\documentclass[11pt,
]{article}

\usepackage{amssymb,amstext,amsthm,latexsym}
\usepackage{amsmath}
\usepackage{amsfonts}

\usepackage{mathrsfs,mathbbol}

\usepackage{mparhack}

\input{q.sty}

\begin{document}

\title{On effective compactness and sigma-compactness}

\author{Vladimir Kanovei\thanks
{
Contact author, 
{\tt kanovei@mccme.ru}}
%
}


\date{\today}
\maketitle


\begin{abstract}
Using the Gandy -- Harrington topology and other methods 
of effective descriptive set theory, we prove several theorems
on compact and \ddd\fsg compact pointsets.
In particular we show that any  
$\is11$ set $A$ of the Baire space $\bn$ either is covered by 
a countable union of compact $\id11$ sets, 
or $A$ contains a subset closed in $\bn$ and homeomorphic to 
$\bn$
(and then $A$ is not covered by a \ddd\fsg compact set, 
of course).
\end{abstract}

\punk{Introduction}

Effective descriptive set theory appeared in the 1950s
as a useful technique of simplification and clarification of 
constructions of classical descriptive set theory 
(see \eg\ \cite{Shon} or \cite{umnKL}). 
Yet it had soon become clear that development of effective 
descriptive set theory leads to results having no direct 
analogies in classical descriptive set theory. 
As an example we recall the following \rit{basis theorem}:
any countable $\id11$ set $A$ of the Baire space 
$\bn=\dN^\dN$ consists of $\id11$ points.
Its remote predecessor in classical descriptive set theory 
is the Luzin -- Novikov theorem on Borel sets with countable
cross-sections.

In this note, methods of effective descriptive set theory 
are applied to the properties of compactness and  
\ddd\fsg compactness of pointsets. 
The following theorem is our main result. 

Recall that $[T]=\ens{x\in\bn}{\kaz m\,(x\res m\in T}$ 
for any tree $T\sq\nse$.

\bte
\lam{mt'}
If\/ $A\sq\bn$ is a\/ $\is11$ set then one and only one of 
the following two claims holds$:$
\ben
\tenu{{\rm(\Roman{enumi})}}
\itla{mt'1}\msur
$A$ is covered by the union\/ $U$ of all sets of the form\/ 
$[T]$, where\/ $T\sq\nse$ is a compact\/ $\id11$ tree 
--- and moreover there is a\/ $\id11$ 
sequence\/ $\sis{T_n}{n\in\dN}$ of compact trees\/ 
$T_n\sq\nse$ such that\/ $A\sq\bigcup_n[T_n]$$;$  

\itla{mt'2}
there is a set\/ $Y\sq A$ homeomorphic to\/ $\bn$ 
and closed in\/ $\bn.$  
\een
\ete

Here conditions \ref{mt'1} and \ref{mt'2} are incompatible: 
if $Y$ is a set as is \ref{mt'2} then $Y$ cannot be 
covered by a \ddd\fsg compact set $U$ as in \ref{mt'1}. 

In parallel to Theorem~\ref{mt'} and using basically the same 
technique, we prove the following similar theorem, which is, 
on the other hand, a  
direct corollary of some well-known results in this field. 

\bte
\lam{mt}
If\/ $A\sq\bn$ is a\/ $\id11$ set then one and only one of 
the following two claims holds$:$
\ben
\tenu{{\rm(\Roman{enumi})}}
\itla{mt1}\msur
$A$ is equal to the union\/ $U$ of all sets of the form\/ $[T]$, 
where\/ $T\sq\nse$ is a compact\/ $\id11$ tree and\/ $[T]\sq A$
--- and moreover there is a\/ $\id11$ 
sequence\/ $\sis{T_n}{n\in\dN}$ of compact trees\/ 
$T_n\sq\nse$ such that\/ $A=\bigcup_n[T_n]$$;$  

\itla{mt2}
there is a set\/ $Y\sq A$ homeomorphic to\/ $\bn$ and relatively 
closed in\/ $A$.
\een
\ete

Conditions \ref{mt1} and \ref{mt2} of the theorem are 
incompatible since $A$ is \ddd\fsg compact provided \ref{mt1} 
holds, so that any relatively closed subset of $A$ is 
\ddd\fsg compact itself, while 
the space $\bn$ is not \ddd\fsg compact, of course.

Theorem~\ref{mt} has strong connections with 4F.18 in \cite{mos} 
which the author of \cite{dst} credits to Louveau. 
It is clear from 4F.18 that if $A$ is a $\id11$ subset 
of $\bn$ and \ddd\fsg compact then it is equal to the 
union of compact $\id11$ sets $A'\sq A$. 
On the other hand, it follows from 4F.14 in \cite{mos} 
that if $A$ is a compact $\id11$ subset of $\bn$ then there 
is a compact $\id11$ tree $T\sq\nse$ such that $A=[T]$. 
To conclude, if $A$ is a \ddd\fsg compact $\id11$ subset 
of $\bn$ then condition \ref{mt1} of Theorem~\ref{mt} is true. 
This allows to derive directly Theorem~\ref{mt}. 
Indeed if $A\sq\bn$ is a $\id11$ set and it does not satisfy
condition \ref{mt1} of Theorem~\ref{mt} then   
the set $A$ is not \ddd\fsg compact by the above, 
and so from the theorem of Hurewicz (see Theorem~\ref{hur})
the set $A$ satisfies \ref{mt2} of Theorem~\ref{mt}.

Nevertheless we present here a new proof of Theorem~\ref{mt}, 
in particular, as a base for the proof of 
\rit{a similar 
but more complicated dichotomy theorem on\/ $\is11$ sets\/} 
(Theorem~\ref{pt}), 
where, unfortunately, the level of effectivity of the 
covering by \ddd\fsg compact sets in \ref{pt1} 
will be less definite.

In addition, we'll prove 
\rit{a generalization of Theorem~\ref{mt'}}
(Theorem~\ref{nt}) which deals, instead of compact sets, with 
closed sets whose trees contain branchings small in the 
sense of a chosen ideal on $\dN$. 

As usual, the theorems remain true in the
relativized form, \ie\ when classes $\id11$ and $\is11$ are  
replaced by $\id11(p)$ and $\is11(p)$, where $p\in\bn$  
is a fixed parameter, with basically the same proofs.

Some well-known classical results related to the theorems 
above are discussed in the last section.

The authors thank anonymous referees for valuable remarks 
and suggestions, including an essential improvement in the 
proof of Theorem~\ref{tks}.

\punk{Preliminaries}
\las{oo}

We use standard notation $\is11\yd\ip11\yd\id11$ for 
effective classes of points and pointsets in $\bn$, as well 
as $\fs11\yd\fp11\yd\fd11$ for corresponding projective classes. 

Let $\nse$ be the set of all finite strings of natural 
numbers, 
including the empty string $\La$. 
If $s,t\in\nse$ then $\lh s$ is the \rit{length} of $s$, 
and $s\su t$ means that $t$ is a 
\rit{proper extension} of $s$. 
If $s\in \nse$ and $n\in\dN$ then $s\we n$ is the string 
obtained by adding $n$ to $s$ as the rightmost term. 
Let, for $s\in\nse,$ 
$$
\ibn s= \ens{x\in\bn}{s\su x}\quad
\text{(\rit{a Baire interval} in $\bn$)}\,. 
$$
If a set $X\sq\bn$ contains at least two elements then there 
is a longest string $s=\stem X$ such that $X\sq\ibn s$.
We put $\diam X=\frac1{1+\stem X}$ in this case, 
and additionally $\diam X=0$ whenever $X$ has at most one 
element.

A set $T\sq\nse$ is a \rit{tree} if 
$s\in T$ holds whenever $s\we n\in T$ for at least one $n$, 
and a \rit{pruned} tree iff $s\in T$ implies $s\we n\in T$ 
for at least one $n$. 
%
Any non-empty tree contains $\La$. 
A string $s\in T$ is a \rit{branching point} of $T$ if 
there are 
$k\ne n$ such that $s\we k\in T$ 
and $s\we n\in T$; let $\bran T$ be the set of all 
branching points of $T$. 
The \rit{branching height} $\bh Ts$ of a string $s\in T$ in a 
tree $T$ is equal to the number of strings 
$t\in\bran T\yt t\su s$.
For instance, if $T=\nse$ then $\bh{\nse}s=\lh s$ 
for any string $s$.
A tree $T$ is \rit{perfect} iff for any $s\in T$ there 
is a string $t\in\bran T$ such that $s\su t$.

A tree $T\sq\nse$ is \rit{compact}, if
it is pruned and has \rit{finite branchings}, that is, 
if $s\in T$ then  
$s\we n\in T$ holds for at most finitely many $n$.
Then
$$
[T]=\ens{x\in\bn}{\kaz m\:(x\res m\in T)} 
$$
is a compact set.
Conversely, if $X\sq\bn$ is a compact set then 
$$
\der X=\ens{x\res n}{x\in X\land n\in\dN}
$$
is a compact tree.
Let $\ct$ be the $\id11$ set of all non-empty 
compact trees.

If $\pX\yi\pY$ are any sets and $P\sq\pX\ti\pY$ then 
$$
\pr P=\ens{x\in\pX}{\sus y\:(\ang{x,y}\in P)}
\quad\text{and}\quad
\seq Px=\ens{y\in\pY}{\ang{x,y}\in P}
$$ 
are, resp., the \rit{projection} of $P$ to $\pX$, and the 
\rit{cross-section} of $P$ defined by $x\in\pX$.
A set $P\sq\pX\ti\pY$ is \rit{uniform} if every 
cross-section $\seq Px$ ($x\in\pX$) contains at most one 
element.
If $P\sq Q\sq \pX\ti\pY$, $P$ is uniform, and 
$\pr P=\pr Q$, then they say that $P$ \rit{uniformizes} $Q$.

\punk{Some facts of effective descriptive set theory}
\las{some}

We'll make use of the following well-known results. 

\bfa
[\rm $\is11$ Separation]
\lam{21}
If\/ $X,Y\sq\bn$ are disjoint\/ $\is11$ sets then there is a\/ 
$\id11$ set\/ $Z\sq\bn$ such that\/ $X\sq Z$ and\/ 
$Y\cap Z=\pu$.\qed
\efa

\bfa
[\rm Kreisel selection, 4B.5 in \cite{mos}]
\lam{22}
If\/ $P\sq\bn\ti\dN$ is a\/ $\ip11$ set and the projection\/ 
$\pr P$ is a\/ $\id11$ set then there is a\/ $\id11$ map\/ 
$f:\pr P\to \dN$ such that\/ $\ang{x,f(x)}\in P$ for all\/ 
$x\in\pr P$.\qed
\efa


\bfa
[\rm 4D.3 in \cite{mos}]
\lam{BQ}
If\/ $P(x,y,z,\dots)$ is a\/ $\ip11$ relation\/  
{\rm(where the domain of each argument can be\/ $\bn$, 
$\cP(\nse)$, the set 
of all compact trees in $\nse$, or any 
other similar Polish space)}
then the relations\/ $\sus x\in\id11\,P(x,y,z,\dots)$ and\/ 
$\sus x\in\id11(y)\,P(x,y,z,\dots)$ 
are\/ $\ip11$, too.\qed
\efa

\bfa
[\rm 4D.14 in \cite{mos}]
\lam{dp}
The set\/ $D=\ens{T\sq\nse}{T\,\text{\rm\ is }\,\id11}$
is\/ $\ip11$.

The set 
$\ens{\ang{p,T}}
{p\in\bn\land T\sq\nse\land T \text{ is } \id11(p)}$ 
is\/ $\ip11$ as well.\qed
\efa

To prove the first claim, note that 
${T\in D}\eqv \sus T'\in\id11\,(T=T')$; 
then the result follows from Fact~\ref{BQ}.

\bfa
[\rm Enumeration of $\id11$ trees, 4D.2 in \cite{mos}]
\lam{23}
There exist\/ $\ip11$ sets\/ $E\sq\dN$ and\/ 
$W\sq\dN\ti\nse$, and a\/ $\is11$ set\/ $W'\sq\dN\ti\nse$ 
such that\/ 
\ben
\renu
\itsep
\itla{-enu1}
$\seq We=\seq {W'}e$ for any\/ $e\in E$\/ 
{\rm(where $\seq {W}e=\ens{s\in\nse}{\ang{e,s}\in W}$)}$;$ 

\itla{-enu2}
a set\/ $T\sq\nse$ is\/ $\id11$ iff there is a number\/ 
$e\in E$ such that\/ $T=\seq We=\seq{W'}e$.\qed
\een
\efa

\vyk{
let $E\yd W\yd W'$ be as in Fact~\ref{23}. 
Then
$$
\bay[b]{rcl}
\hspace*{0ex}
T\in D 
\hspace{-1ex}&\eqv&\hspace{-1ex} 
\sus e\in E\:(T=\seq We)\\[\dxii]
\hspace{-1ex}&\eqv&\hspace{-1ex} 
\sus e\in E\; 
\kaz s\in\nse\:(s\in\seq{W'}e\imp s\in T\imp s\in\seq We). 
\eay
$$
}

\bfa
[4F.17 in \cite{mos}]
\lam{25}
If\/ $P\sq\bn\ti\bn$ is a\/ $\id11$ set and every 
cross-section\/ $\seq Px$ ($x\in\bn$) is at most countable then
$\pr P$ is a\/ $\id11$ set, and\/
$P$ is a countable union of uniform\/ $\id11$ sets each of which 
uniformizes\/ $P$.\qed
\efa

\bfa
[4F.14 in \cite{mos}]
\lam{ks=t}
If\/ $F\sq\bn$ is a closed\/ $\id11$ set and\/  
$X\sq F$ is a compact\/ $\is11$ set then there is a
compact\/ $\id11$ tree\/ $T\sq\nse$ such that\/ 
$X\sq[T]\sq F$. 
In particular, in the case\/ $X=F$, any compact\/ 
$\id11$ set\/ $X\sq\bn$ has the form\/ $X=[T]$ for 
some compact\/ $\id11$ tree\/ $T\sq\nse$. 
\efa

\vyk{
\bpf
(Added for the sake of the reader's convenience.) 
Then the set
$$
H=\ens{\ang{t,n}}{t\in\nse\land n\in\dN\land
\kaz k\:(t\we k\in T(Y)\imp k \le n)}
$$
belongs to $\ip11$ 
(as the $\is11$ set $T(Y)$ stands to the left of an implication), 
and $\dom H=\nse$. 
Fact~\ref{22} implies that there is a $\id11$ map 
$f:\nse\to\dN$ such that $\ang{t,f(t)}\in H$ for all $t\in\nse$. 
We have $y(n)\le f(y\res n)$ for all $y\in Y$ and $n$
by the definition of $H$, and hence $Y\sq[T']$, 
where $T'$ consists of all strings $t\in\nse$ satisfying 
$t(n)\le f(t\res n)$ for every $n<\lh t$.
Yet $T'$ is a compact tree, and clearly $T'$ is $\id11$ 
by the choice of $f$. 
Therefore $[T']\sq U$, and this is a contradiction.
\epF{Lemma}
}

\bfa
[4F.11 in \cite{mos}]
\lam{DinD}
Any compact\/ $\id11$ set\/ $\pu\ne A\sq\bn$ contains a\/ 
$\id11$ element\/ $x\in A$.\qed
\efa   

There is a useful uniform version of Fact~\ref{23}.

\bfa
[\rm Uniform enumeration]
\lam{23+}
There exist\/ $\ip11$ sets\/ $\bE\sq\bn\ti\dN$ and\/ 
$\bW\sq\bn\ti\dN\ti\nse$, 
and a\/ $\is11$ set\/ $\bW'\sq\bn\ti\dN\ti\nse$ 
such that\/ 
\ben
\renu
\itsep
\itla{+enu1}
$\sek \bW xe=\sek{\bW'}xe$ for any\/ $\ang{x,e}\in \bE$\/ \ 
{\rm(where $\sek \bW xe=\ens{s\in\nse}{\ang{x,e,s}\in \bW}$)}$;$ 

\itla{+enu2}
if\/ $x\in\bn$ then
a set\/ $T\sq\nse$ is\/ $\id11(x)$ iff there is a number\/ 
$e\in E$ such that\/ $T=\sek \bW xe=\sek {\bW'} xe$.\qed
\een
\efa

This result allows us to prove the following generalization 
of Fact~\ref{22}, also well-known in effective descriptive 
set theory.

\bfa
[\rm 4D.6 in \cite{mos}]
\lam{22+}
Suppose that\/ $Q\sq\bn\ti\cP(\nse)$ is\/ $\ip11$, 
the projection\/ $\pr Q$ onto\/ $\bn$ is\/ $\id11$, 
and for each\/ $x\in\pr Q$ there exists a set\/ $T\in\id11(x)$ 
such that\/ $\ang{x,T}\in Q$. 
Then there is a\/ $\id11$ map\/ 
$\tau:\pr Q\to \cP(\nse)$ such that\/ $\ang{x,\tau(x)}\in Q$ 
for all\/ $x\in\pr Q$. 
\efa
\bpf
Making use of sets $\bE,\bW,\bW'$ as in Fact~\ref{23+}, we let 
$$
P=\ens{\ang{x,e}\in\bE}{\ang{x,\sek \bW xe}\in Q}.
$$
Immediately the set $P$ is $\ip11$ and $\pr P=\pr Q$ is a 
$\id11$ subset of $\bn$. 
By Fact~\ref{22}, there is a $\id11$ map $f:\pr P\to \dN$ 
such that $\ang{x,f(x)}\in P$ for all\/ $x\in\pr P$. 
It remains to define $\tau(x)=\sek \bW x{f(x)}$ for all 
$x\in\pr Q$; to prove that $\tau$ is $\id11$  
use both $\bW$ and $\bW'$.
\epf

Facts 
\ref{21}, \ref{22}, \ref{BQ}, \ref{23}, 
\ref{dp} (the first claim), 
\ref{25}, \ref{ks=t}, \ref{23+} , \ref{22+} 
remain true for relativized lightface classes  
$\is11(p)\yd\ip11(p)\yd\id11(p)$, where $p\in\bn$ is 
an arbitrary fixed parameter.
Therefore Facts \ref{21}, \ref{22}, \ref{25} also 
hold with lightface classes replaced by boldface 
projective classes $\fs11\yd\fp11\yd\fd11$.

\punk{The Gandy -- Harrington topology}
\las{gaha1}

The \rit{\gh\ topology} on the Baire space $\bn$ consists of  
all unions of $\is11$ sets $S\sq\bn$.
This topology includes the Polish topology on $\bn$ but 
is not Polish. 
Nevertheless the \gh\ topology satisfies a condition  
typical for Polish spaces.

\bdf
\lam{genb'}
Let $\cF$ be any family of sets, \eg\ sets in a given 
background space $\dX$.
A set $\cD\sq\cF$ is \rit{open dense}   
iff \,
$\kaz F\in\cF\:\sus D\in\cD\:(D\sq F)$, \, and %
$$
\kaz F\in\cF\:\kaz D\in\cD\:(F\sq D\imp F\in\cD)\,.
$$
Sets $\cD$ satisfying only the first requirement  
are called \rit{dense}.
If $\cD\sq\cF$ is dense then the set 
$\cD'=\ens{F\in\cF}{\sus D\in\cD\:(F\sq D)}$ is open dense. 
The notions of \rit{open} and \rit{dense} are related to 
a certain topology which we'll not discuss, but not necessarily 
with the topology of the background space $\dX$. 

A \rit{Polish net} for $\cF$ is any   
collection $\ens{\cD_n}{n\in\dN}$ of open dense sets 
$\cD_n\sq\cF$ such that we have $\bigcap_nF_n\ne\pu$ 
for every sequence of sets 
$F_n\in\cD_n$ satisfying the finite intersection property  
(\ie\ $\bigcap_{k\le n}F_k\ne\pu$ for all $n$). 
\edf

For instance the family of all non-empty closed sets of a 
complete metric space $\pX$ admits a Polish net: 
let $\cD_n$ contain all closed sets of diameter 
$\le n\obr$ in $\pX$.
The next theorem is less elementary. 
This theorem and the following corollary are well-known, see \eg\  
\cite{hkl,h-ban,umnS,k-ams}. 

\bte
\lam{s11p}
The collection\/ $\dP$ of all non-empty\/ $\is11$ sets  
in\/ $\bn$ admits a Polish net.\qed
\ete

\vyk{
\bcor
\lam{s11pc}
The space\/ $\bn$ with the \gh\ topology satisfies the  
Baire theorem, \ie, all comeager sets are dense.
\qed
\ecor
}

\vyk{
\punk{On compact $\id11$ sets}
\las{gahaK}

If a tree $T\sq\nse$ is $\id11$ then the set 
$[T]=\ens{a\in\bn}{\kaz m\:(a\res m\in T)}$ is $\id11$ 
as well since $x\in[T]\eqv \kaz m\:(x\res m\in T)$. 

\ble
\lam{ks=t}
If\/ $F\sq\bn$ is a closed\/ $\id11$ set and\/  
$A\sq F$ is a compact\/ $\is11$ set then there is a
compact\/ $\id11$ tree\/ $T\sq\nse$ such that\/ 
$A\sq[T]\sq F$. 

In particular, in the case\/ $A=F$, any compact\/ 
$\id11$ set\/ $A\sq\bn$ has the form\/ $A=[T]$ for 
some compact\/ $\id11$ tree\/ $T\sq\nse$.
\ele

\bpf
\rit{Part 1}.
We claim that there is a $\id11$ tree $S$ 
(not necessarily a compact tree) such that 
$F=[S]$. 
Indeed the complementary $\id11$ set $G=\bn\bez F$ is open. 
Therefore the set
$$
P=\ens{\ang{x,t}\in \bn\ti\nse}{x\in\ibn t\sq G}
$$
satisfies $\dom P=G$.
Moreover $P$ is $\ip11$ since the relation 
$\ibn t\sq G$ can be expressed by a $\ip11$ formula 
$\kaz x\:(t\su x\imp x\in G)$.  
Fact~\ref{22} 
yields a $\id11$ map $f:G\to\nse$ such that 
$x\in\ibn{f(x)} \sq G$ for each $x\in G$. 
Then the set 
$$
U=\ens{f(x)}{x\in G}=\ens{t\in\nse}{\sus x\in G\:(f(x)=t)}
$$
belongs to $\is11$ and satisfies 
$G=\bigcup_{t\in U}\ibn t$. 
However
$$
V=\ens{t\in\nse}{\ibn t\sq G}=
\ens{t\in\nse}{\kaz x\:(x\in\ibn t\imp x\in G)}
$$
is a $\ip11$ set, $U\sq V$, and 
$G=\bigcup_{t\in U}\ibn t$. 
It follows from Fact~\ref{21} that there is a 
$\id11$ set $W$ such that $U\sq W\sq V$ --- and then 
$G=\bigcup_{t\in W}\ibn t$. 

Now put $S=\ens{s\in\nse}{\kaz t\:(t\in W\imp t\not\sq s)}$. 
Easily $S$ is a $\id11$ tree 
(possibly with endpoints), and $[S]=F$. 
This completes part 1 of the proof.\vom

\rit{Part 2}. 
The set $P$ of all pairs $\ang{s,u}$ such that $s\in\nse$, 
$u\sq\dN$ is finite and non-empty, and 
$$
\kaz x\in\bn\:
\big((x\in A\land
s\su x)\imp \sus k\in u\:
(s\we k\su x)\land \kaz k\in u\:(s\we k\in S)
\big),
$$
is $\ip11$ in $\dN\ti\pwf\dN$, where the second factor 
(= the set of all finite $u\sq\dN$)
is identified with $\dN$ by means of a recursive  
bijection.
Moreover $\dom P=\nse$. 
(Let $s\in\nse$, $n=\lh s$. 
If there is no $x\in A$ with 
$s\su x$, then $\ang{s,u}\in P$ for every finite $u$. 
If there is $x\in A$ with $s\su x$ then the set 
$u=\ens{x(n)}{s\su x\in A}$ is finite by the compactness, 
and hence $\ang{s,u}\in P$.) 
Therefore, as above, there is a $\id11$ map $f:\nse\to\pwf\dN$ 
such that $\ang{s,f(s)}\in P$ for each $s\in\nse$. 

The tree $T=\ens{s\in\nse}{\kaz n<\lh s\:(s(n)\in f(s\res n))}$ 
has finite branchings as all values of $f$ are finite sets. 
To show that $T$ has no endpoints,  
let $s\in T$. 
Then $s(n)\in f(s\res n)$ for all $n<\lh s$,  
and $f(s)$ is non-empty. 
Pick any $k\in f(s)$; 
the extended string $s'=s\we k$ belongs to $T$. 
Thus $T$ is a compact $\id11$ tree. 

To prove that $A\sq[T]$ assume that $x\in A$. 
As $\La=x\res0$ (the empty string) obviously belongs to $T$, 
it suffices to prove that if $m<\om$ and $s=x\res m\in T$ 
then the extended string $t=x\res{(m+1)}=s\we x(m)$ belongs to 
$T$ as well, or, equivalently, that $x(m)\in f(s)=f(x\res m)$ 
--- but this holds because $\ang{s,f(s)}\in P$ by the choice 
of $f$.
Finally to prove that $[T]\sq F$ note that by definition 
$T\sq S$ and $[S]=F$.
\epf
}

\vyk{
\bcor
\lam{scomK}
Any non-empty compact\/ $\id11$ set\/ $A\sq\bn$ contains a\/ 
$\id11$ element\/ $x\in A$.
\ecor
\bpf
By Lemma~\ref{ks=t}, there is a compact $\id11$ tree 
$T\sq\nse$ such that $A=[T]$.
Let $x$ the lexicographically leftmost branch in $T$, 
that is,
\pagebreak[1] 
$$
x(n)=\tmin\ens{s(n)}
{s\in T\land n<\lh s\land s\res n=\ang{x(0),\dots,x(n-1)}}, 
\quad\kaz n.
\eqno\qed
$$
\ePf
}

\punk{The proof of Theorem~\ref{mt'}}
\las{d2}

Recall that $\ct$ is the set of all compact trees 
$\pu\ne T\sq\nse$; $\ct$ is $\id11$, of course.
Let $U$ be the set as in \ref{mt'1} of the theorem. 
We claim that $U$ is $\ip11$. 
Indeed,  by definition
$$
x\in U
\leqv
\sus T\in\id11\:
(T\in\ct\land x\in [T])\,,
$$
and the result follows from Fact~\ref{BQ}. 

It follows that the difference $A'=A\dif U$ is a $\is11$ set.

\ble
\lam{tkm*}
Under the conditions of Theorem~\ref{mt'}, 
if\/ $Y\sq A'$ is a non-empty\/ $\is11$ set then 
its topological closure\/ $\clo Y$ in\/ $\bn$ 
is not compact, \ie, the tree\/  
$\der Y =\ens{y\res n}{y\in Y\land n\in\dN}$ 
has at least one infinite branching.
\ele
\bpf
Suppose otherwise: $\clo Y$ is compact. 
Then by Fact~\ref{ks=t} (with $F=\bn$) 
there is a compact $\id11$ tree $T$ such that 
$\clo Y\sq[T]$.
Therefore $Y\sq\clo Y\sq [T]\sq U$, and this contradicts 
to the assumption $\pu\ne Y\sq A'$.
\epF{Lemma}

\rit{Case 1}: the set
$A'=A\dif U$ is non-empty. 
We assert that then there is a system of non-empty 
$\is11$ sets $Y_s\sq A'$ satisfying the following conditions 
\ben
\itsep
\tenu{(\arabic{enumi})}
\itla{gan1}
if $s\in\nse$  and $i\in\dN$ then 
$Y_{s\we i}\sq Y_s$;

\itla{gan2} 
\vyk{
the diameter\footnote
{The diameter of a set $Y\sq\bn$ is equal to $0$ whenever $Y$ 
contains at most one point, and is equal to $\frac1{n+1}$ 
whenever $Y$ contains at least two points and $n$ is the 
largest number such that $x\res n=y\res n$ for all 
$x,y\in Y$.} 
of $Y_s$ is not greater than 
}%
$\diam{Y_s}\le2^{-\lh s}$; 

\itla{gan3}
if $s\in\nse$ and $k\ne n$ then
$Y_{s\we k}\cap Y_{s\we n}=\pu$, and moreover, sets 
$Y_{s\we k}$ are covered by pairwise disjoint (clopen) 
Baire intervals $J_{s\we k}$;

\itla{gan4} 
$Y_s\in\cD_{\lh s}$, where by Theorem~\ref{s11p} 
$\ens{\cD_n}{n\in\dN}$ is a Polish net for the family  
$\dP$ of all non-empty $\is11$ sets $Y\sq\bn$;

\itla{han5} 
if $s\in\nse$ and $x_k\in Y_{s\we k}$ for all $k\in\dN$ 
then the sequence of points $x_k$ does not have 
convergent subsequences in $\bn$.
\een

If such a construction is accomplished then \ref{gan4} implies 
that $\bigcap_mY_{a\res m}\ne\pu$ for each $a\in\bn$. 
On the other hand by \ref{gan2} every such an intersection 
contains a single point, which we denote by $f(a)$, and the  
map $f:\bn\na Y=\ran f=\ens{f(a)}{a\in\bn}$ is   
a homeomorphism by clear reasons. 

Prove that $Y$ is closed in $\bn$. 
Consider an arbitrary sequence of points $a_n\in\bn$ such 
that the corresponding sequence of points $y_n=f(a_n)\in Y$ 
converges to a point $y\in\bn$; we have to prove that 
$y\in Y$.
If the sequence $\sis{a_n}{n\in\dN}$ contains a subsequence 
of points $b_k=a_{n(k)}$ convergent to some 
$b\in\bn$ then quite obviously the sequence of points 
$z_k=f(b_k)$ (a subsequence of $\sis{y_n}{n\in\dN}$)
converges to $z=f(b)\in Y$, as required. 
Thus suppose that the sequence $\sis{a_n}{n\in\dN}$ 
has no convergent subsequences. 
Then it cannot be covered by a compact set, and it easily 
follows that there is a string $s\in\nse$, an infinite 
set $K\sq\dN$, and for each $k\in K$ --- a number 
$n(k)$ such that $s\we k\su a_{n(k)}$. 
But then $y_{n(k)}\in Y_{s\we k}$ by construction. 
Therefore the subsequence $\sis{y_{n(k)}}{k\in\dN}$ 
diverges by \ref{han5}, which is a contradiction.

Thus $Y$ is closed,
and hence we have \ref{mt'2} of Theorem~\ref{mt'}.  

As for the construction of sets $Y_s$, 
if a $\is11$ set $Y_s\sq A'$ is defined then 
by Lemma~\ref{tkm*} there is a string $t\in T(Y_s)$ such 
that $t\we k\in T(Y_s)$ for all $k$ in an infinite set 
$K_s\sq\dN$. 
This allows us to define a sequence of pairwise different 
points $y_k\in Y_s$ ($k\in\dN$) having no convergent 
subsequences. 
We cover these points by Baire intervals $U_k$ small enough 
for \ref{han5} to be true for the $\is11$ sets  
$Y_{s\we i}=Y_s\cap U_i$, and then shrink these sets if 
necessary to satisfy \ref{gan2} and \ref{gan4}.\vom

\rit{Case 2}: $U=\pu$, that is, $A\sq U$. 
Recall that $\ct$ is the $\id11$ set of all compact trees 
$T\sq\nse$. 
The sets
$$
\bay{cclcr}
Q&=&\ens{\ang{x,T}}
{x\in\bn\land T\in \ct\cap\id11\land x\in[T]}\,,
&&\text{and\,}\\[1ex]
Z&=&\ens{x\in\bn}
{\sus T\in\id11\,(T\in\ct\land x\in[T])}&=&\pr Q
\eay
$$
are $\ip11$ by Facts~\ref{BQ} and \ref{dp}. 
Moreover, $A\sq U$ implies  
$A\sq Z$, and hence by Fact \ref{21} there is a $\id11$ set 
$X$ such that $A\sq X\sq Z$. 
Then $P=\ens{\ang{x,n}\in Q}{x\in X}$ is still a $\ip11$ set, 
and $\pr P=X$ is a $\id11$ set. 
Therefore by Fact~\ref{22+} there is a $\id11$ function 
$\tau:X\to\ct$ such that $\ang{x,\tau(x)}\in Q$ for all $x\in X$. 

Note that $\tau(x)\in\ct\cap\id11$ and $x\in [\tau(x)]$ for all 
$x\in A$ by the construction. 
Thus the full image $R=\ens{\tau(x)}{x\in A}$ is a $\is11$ 
subset of the $\ip11$ set $\ct\cap\id11$, 
and hence there is a $\id11$ set $D$ such that 
$R\sq D\sq \ct\cap\id11$. 
But countable $\id11$ sets are known to admit a $\id11$ 
enumeration, so there is a $\id11$ map $\delta:\dN\onto D$.
Now let $T_n=\delta(n)$ for for all $n$.

\qeDD{Theorem~\ref{mt'}}

\punk{The proof of Theorem~\ref{mt}}
\las{gahaSK}

By Theorem~\ref{mt'}, we can \noo\ assume 
that $A$ is covered by a \ddd\fsg compact set, and hence 
if $F\sq A$ is a closed set then $F$ is \ddd\fsg compact.
Further, the set $U$ in \ref{mt1} of Theorem~\ref{mt} 
(the union of all sets $[T]\sq A$, 
where $T$ is a compact $\id11$ tree) 
is $\ip11$.
Indeed,
$$
x\in U
\leqv
\sus T\in\id11\:
(T\,\text{ is a compact tree and }\,x\in [T]\sq A)\,,
$$
and the result follows from Fact~\ref{BQ} since 
the property of \lap{being a compact tree} 
can be straightforwardly expressed by an arithmetic formula,  
while $[T]\sq A$ can be expressed by a $\ip11$ formula.

We conclude that $A'=A\bez U$ is $\is11$. 

\ble
\lam{tkm-l}
If\/ $F\sq A'$ is a non-empty\/ $\is11$ set then\/ 
$\clo F\not\sq A$. 
\ele
\bpf
We first prove that if $X\sq A$ is a compact $\is11$ set 
then $A'\cap X=\pu$.
Suppose towards the contrary that $A'\cap X$ is non-empty. 
We are going to find a closed $\id11$ set $F$ satisfying 
$X\sq F\sq A$ --- this would imply $X\sq U$ by 
Fact~\ref{ks=t}, which is a contradiction. 

Since the complementary $\ip11$ set $C=\bn\bez X$ is open, 
the set 
$$
H=\ens{\ang{x,s}}
{s\in\nse\land x\in C\cap\ibn s\land \ibn s\cap X=\pu}
$$ 
is $\ip11$ and $\pr H=C$. 
Thus the $\id11$ set $D=\bn\bez A$ is included 
in $\pr H$. 
By Fact~\ref{22},
there is a $\id11$ map $\nu:D\to\nse$ such that 
$x\in D\imp \ang{x,\nu(x)}\in H$, or equivalently, 
$x\in \ibn{\nu(x)}\sq C$ for all $x\in D$. 
Then the set $\Sg=\ran\nu=\ens{\nu(x)}{x\in D}\sq\nse$ 
is $\is11$ and $D\sq \bigcup_{s\in\Sg}\ibn s\sq C$.

But $\Pi=\ens{s\in\nse}{\ibn s\sq C}$ is a $\ip11$ set
and $\Sg\sq \Pi$.
It follows that there exists a 
$\id11$ set $\Da$ such that $\Sg\sq\Da\sq\Pi$. 
Then still $D\sq \bigcup_{s\in\Da}\ibn s\sq C$,
and hence the closed set 
$F=\bn\bez \bigcup_{s\in\Da}\ibn s$ satisfies $X\sq F\sq A$. 
But $x\in F$ is equivalent to $\kaz s\:(s\in\Da\imp x\nin \ibn s)$, 
thus $F$ is $\id11$, as required.

Now suppose towards the contrary that 
$\pu\ne F\sq A'$ is a $\is11$ set but $\clo F\sq A$. 
By the \noo\ assumption above,  
$\clo F=\bigcup_nF_n$ is \ddd\fsg compact, 
where all $F_n$ are compact. 
There is a Baire interval $\ibn s$ such that the set 
$X=\ibn s\cap \clo F$ is non-empty and $X\sq F_n$ for some $n$. 
Thus $X\sq A$ is a non-empty compact $\is11$ set, hence 
$X\cap A'=\pu$ by the first part of the proof. 
In other words, $\ibn s\cap \clo F\cap A'=\pu$.
It follows that $\ibn s\cap F=\pu$ (because $F\sq A'$), 
contrary to $X=\ibn s\cap \clo F\ne \pu$. 
\epF{Lemma}

\rit{Case 1}: 
the $\is11$ set $A'\sq A$ is non-empty. 
To get a set $Y\sq A'$, \rit{relatively} 
closed in $A$ and homeomorphic to $\bn$, as 
in \ref{mt2} of the theorem, we'll define  
a system of non-empty $\is11$ sets $Y_s\sq A'$ satisfying 
conditions \ref{gan1}, \ref{gan2}, \ref{gan3}, \ref{gan4} of 
Section~\ref{d2}, along with the next requirement 
instead of \ref{han5}:
\ben
\itsep
\tenu{$(\arabic{enumi}')$}
\atc\atc\atc\atc
\itla{gan5} 
if $s\in\nse$ then there is a point $y_s\in \clo{Y_s}\bez A$ 
such that any sequence of points 
$x_k\in Y_{s\we k}$ ($k\in\dN$) converges to $y_s$.
\een

If we have defined such a system of sets, then the associated 
map $f:\bn\to A'$ is $1-1$ and is a homeomorphism 
from $\bn$ onto its full image 
$Y=\ran f =\ens{f(a)}{a\in\bn}\sq A'$,
as in the proof of Theorem~\ref{mt'}.  

Let's prove that $Y$ is relatively closed in $A$. 
Consider a sequence of points $a_n\in\bn$ such 
that the corresponding sequence of $y_n=f(a_n)\in Y$ 
converges to a point $y\in\bn$; 
we have to prove that $y\in Y$ or $y\nin A$.
If the sequence $\sis{a_n}{}$ contains a subsequence 
convergent to $b\in\bn$ then, as in the proof of 
Theorem~\ref{mt'}, 
$\sis{y_n}{}$ converges to $f(b)\in Y$. 
If the sequence $\sis{a_n}{}$ 
has no convergent subsequences, 
then there exist a string $s\in\nse$, 
an infinite set $K\sq\dN$, 
and for each $k\in K$ --- a number $n(k)$, 
such that $s\we k\su a_{n(k)}$. 
But then $y_{n(k)}\in Y_{s\we k}$ by construction. 
Therefore the subsequence $\sis{y_{n(k)}}{k\in\dN}$ 
converges to a point $y_s\nin A$ by \ref{han5}, 
as required.

Finally on the construction of sets $Y_s$. 

Suppose that a $\is11$ set $\pu\ne Y_s\sq A'$ is defined. 
Then its closure $\clo{Y_s}$ is a $\is11$ set, too, 
therefore $\clo{Y_s}\not\sq A$ by Lemma~\ref{tkm-l}.
There is a sequence of 
pairwise different points $x_n\in Y_s$ which converges to 
a point $y_s\in \clo{Y_s}\bez A$. 
Let $U_n$ be a neighbourhood of $x_n$ (a Baire interval) 
of diameter less than $\frac13$ of the least distance 
from $x_n$ to the points $x_k\yt k\ne n$.
Put $Y_{s\we n}=Y_s\cap U_n$, and shrink the sets 
$Y_{s\we n}$ so that they satisfy \ref{gan2} 
and \ref{gan4}.\vom

\rit{Case 2}: $A'=\pu$, that is, $A=U$. 
This implies \ref{mt1} of the theorem, exactly as 
in the proof of Theorem~\ref{mt'} above.\vom
 
\qeDD{Theorem~\ref{mt}}

\punk{A generalization of Theorem~\ref{mt'}}
\las{gen}

Let $\cI\sq\pws\dN$ be an ideal on $\dN$. 
A tree $T\sq\nse$ is:
\bde
\item[\it\ddi small,]
if for any $s\in T$ the set 
$\suc Ts=\ens{n}{s\we n\in T}$ belongs to $\cI$;

\item[\it\ddi positive,]
if 1) it is perfect, and 
2) if $s\in\bran T$ then the set 
$\suc Ts$ does \rit{not} belong to $\cI$.
\ede
Accordingly, a set $X\sq\bn$ is:
\bde
\item[\it\ddi small,]
if $\der X=\ens{x\res n}{n\in\dN\land x\in X}$ is an 
\ddi small tree;

\item[\it\dsi small,]
if it is a countable union of \ddi small sets;

\item[\it\ddi positive,]
if it contains a subset of the form $[T]$, 
where $T\sq\nse$ is an \ddi positive tree.
\ede
For instance, if $\cI=\fin$ is the Frechet ideal of all 
finite sets $x\sq\dN$ then \ddi small trees and sets are 
exactly compact trees, resp., sets, 
\dsi small sets are \ddd\fsg compact sets,
while \ddi positive trees are perfect trees 
with infinite branchings.
Moreover if $T$ is such a \ddd\fin positive tree then the 
set $[T]$ is closed and homeomorphic to $\bn$, hence, 
non-\ddd\fsg compact.
Thus condition \ref{mt'2} of Theorem~\ref{mt'} can be 
reformulated as follows:
\rit{$A$ is a\/ \ddd\fin positive set.}

Here we prove the following theorem 
(compare with Theorem~\ref{mt'}).

\bte
\lam{nt}
Let\/ $\cI$ be a\/ $\ip11$ ideal on\/ $\dN$.
If\/ $A\sq\bn$ is a\/ $\is11$ set then one and only one of 
the following two claims holds$:$
\ben
\tenu{{\rm(\Roman{enumi})}}
\itla{nt1}\msur
$A$ is\/ \dsi small$;$ 

\itla{nt2}
$A$ is an\/ \ddi positive set.
\een
\ete

Condition \ref{nt1} of this theorem is notably weaker than 
a true generalization of Theorem~\ref{mt'} would require:
\rit{$A$ is covered by the union of all sets\/ $[T]$, 
where\/ $T\sq\nse$ is an\/ \ddi small\/ $\id11$ tree.}
Unfortunately such a stronger version is not accessible 
so far. 
The key element in the proof of Theorem~\ref{mt'}, which 
allows to strengthen \ref{mt'1} from $\is11$ to $\id11$, 
is Lemma~\ref{tkm*} based on Fact~\ref{ks=t}.
We don't know whether the latter is true in the context 
of Theorem~\ref{nt}, \eg, at least in the form: 
\rit{any \ddi small\/ $\is11$ set is covered by a\/ 
\ddi small\/ $\id11$ set}. 
It would be sufficient to assume that $\cI$ satisfies the 
following: 
\rit{
if\/ $p\in\bn$ and\/ $x\in\cI$ is a\/ $\is11(p)$ set then 
there is a\/ $\id11(p)$ set\/ $y\in\cI$ such that\/ $x\sq y$.}

\bpf
As covering of small $\is11$ sets by small $\id11$ sets is 
not available, we'll follow a line of arguments which differ 
from the proof of Theorem~\ref{mt'} above.
First of all, $A=\pr P=\ens{x\in\bn}{\sus y\:P(x,y)}$, 
where $P\sq\bn\ti\bn$ is a $\ip01$ set. 
Consider the tree
$$
S=\ens{\ang{x\res n,y\res n}}
{n\in\dN\land \ang{x,y}\in P} \sq\nse\ti\nse,
$$
so that 
$P=[S]=\ens{\ang{x,y}\in\bn^2}
{\kaz n\:(\ang{x\res n,y\res n}\in S)}$.
If $u,v\in\nse$ then let 
$P_{uv}=\ens{\ang{x,y}\in P}{u\su x\land v\su y}$ and 
$A_{uv}=\pr P_{uv}$; 
thus, in particular, $P_{\La\La}=P$ and $A_{\La\La}=A$. 
If the subtree 
$$
S'=\ens{\ang{u,v}\in S}
{A_{uv}\text{ is not \dsi small}}
$$ 
of $S$ is empty then $A=A_{\La\La}$ is \dsi small, 
getting \ref{nt1} of the theorem. 
Therefore we assume that $S'\ne\pu$, and the goal is 
to get \ref{nt2} of the theorem.  

Note that $P_{uv}=\bigcup_{k,n}P_{u\we k,v\we n}$, and hence 
the tree $S'$ has no maximal nodes: if $\ang{u,v}\in S'$ then 
$\ang{u\we k,v\we n}\in S'$ for some $k\yi n$.
We consider the corresponding closed set  
$$
P'=[S']=\ens{\ang{x,y}\in\bn^2}
{\kaz n\:(\ang{x\res n,y\res n}\in S')}
$$
and the $\fs11$ set $A'=\pr{P'}$. 
If $\ang{u,v}\in S'$ then let 
$$
P'_{uv} =\ens{\ang{x,y}\in P'}
{u\su x\land v\su y}
\quad\text{and}\quad
A'_{uv}=\pr P'_{uv}\,,
$$
so that $A'_{uv}$ is a non-empty $\fs11$ subset of $A'$, 
not \dsi small by the definition of $S'$. 
The next lemma is quite obvious.

\ble
\lam{tt'}
If\/ $\ang{u,v}\in S'$, $u'\in\nse$, $u\su u'$, and\/ 
$A'_{uv}\cap\ibn{u'}\ne\pu$ then there is a string\/ 
$v'\in\nse$ such that\/ $v\su v'$ and\/ 
$\ang{u',v'}\in S'$.\qed
\ele

We are going to define a pruned tree $T\sq\nse$ and a 
string $v(t)\in\nse$ for all $t\in T$, such that  
%
\ben
\itsep
\tenu{(\arabic{enumi})}
\itla{lx1}
if $t\in T$ then $\ang{t,v(t)}\in S'$; 

\itla{lx2}
if $s,t\in T$ and $s\sq t$ then 
$v(s)\sq v(t)$;

\itla{lx3}
if $s\in T$ then there exists a string $t\in T$ such that
$s\su t$ and the set $\ens{k}{t\we k\in T}$ does not 
belong to $\cI$.
\een
If such construction is accomplished then $T$ is an 
\ddi positive tree by \ref{lx3}, 
and on the other hand $[T]\sq A'\sq A$, 
so that \ref{nt2} of the theorem holds. 

Thus it remains to carry out the construction. 

To begin with we define $\La\in T$, of course, and let 
$v(\La)=\La$.

Suppose that $t\in T$, so that $\ang{t,v(t)}\in S'$ and  
the set $A'_{t,v(t)}$ is not \dsi small,
in particular, not \ddi small, hence the tree 
$\der{A'_{t,v(t)}}$ is not \ddi small. 
We conclude that there is a string $s\in\nse$ such that 
$t\sq s$ and the set 
$K=\ens{k}{\sus a\in A'_{t,v(t)}\,(s\we k\su a)}$
does not belong to $\cI$. 

We let every string $t'$ with $t\su t'\sq s$ belong to $T$, 
and choose $v(t')$ for any such $t'$ so that \ref{lx1} and 
\ref{lx2} hold, using Lemma~\ref{tt'}. 
Then let every string $s\we k\yt k\in K$, belong to $T$, 
and let $v(s\we k)=v$, where $v$ is any string such that 
$v(s)\sq v$ and $\ang{s\we k,v}\in S'$.
(The existence of at least one such string $v$ follows from 
Lemma~\ref{tt'}.)\vom

\epF{Theorem~\ref{nt}}


\punk{Theorem~\ref{mt} for $\is11$ sets}
\las{mtS}

There is a difference between Theorem~\ref{mt'} and 
Theorem~\ref{mt}:
the first theorem deals with $\is11$ sets while the other 
one --- with $\id11$ sets. 
We don't know whether Theorem~\ref{mt} holds for all $\is11$ 
sets, but it is quite clear where the proof in 
Section~\ref{gahaSK} fails.  
Indeed if $A$ is a $\is11$ set then $A'$ turns out to be 
a set in $\fs11$ and $\is12$, but not $\is11$, 
so the rest of the proof just does not work. 
Nevertheless we can prove the following essentially weaker 
result.

\bte
\lam{pt}
If\/ $A\sq\bn$ is a\/ $\is11$ set then one and only one of 
the following two claims holds$:$
\ben
\tenu{{\rm(\Roman{enumi})}}
\itla{pt1} 
there exist\/$:$ 
a countable ordinal\/ $\la$ and an effectively defined 
sequence\/ $\sis{T^\al}{\al<\la}$ of compact\/ 
$\id13$ trees\/ $T^\al\sq\nse$ 
such that\/ $A=\bigcup_{\al<\la}[T^\al]$ --- 
{\rm then $A$ is \ddd\fsg compact, of course}$;$ 

\itla{pt2}
there is a set\/ $Y\sq A$ homeomorphic to\/ $\bn$ and relatively 
closed in\/ $A$.
\een
\ete

We'll not try to estimate the level and character of the 
effectivity condition in \ref{pt1}, since we don't think that 
our construction gives a result even close to optimal. 
But it will be quite clear from the construction that it 
is absolute for all transitive models 
containing the true $\omi$, and lies 
within the projective hierarchy and probably within $\id13$.
It is still an interesting {\bf problem} to prove 
Theorem~\ref{mt}, as it stands, for $\is11$ sets.   

\bpf
By Theorem~\ref{mt'}, we can \noo\ assume 
that $A$ is covered by a \ddd\fsg compact set, and hence 
if $F\sq A$ is a closed set then $F$ is \ddd\fsg compact.
Let $P\sq\bn\ti\bn$ be a $\ip01$ set such that 
$A=\pr P$, and
$$
S=\ens{\ang{x\res n,y\res n}}
{n\in\dN\land \ang{x,y}\in P} \sq\nse\ti\nse,
$$
so that $P=[S]$.
A decreasing sequence of derived trees 
$\sa S\al,\;\al\in\Ord$, is 
defined by induction so that 
$\sa S0=S$, if $\la$ is limit then 
$\sa S\la=\bigcap_{\al<\la}\sa S\al$, and 
for any $\al$:  
\ben
\Aenu
\itla{+1A}
we let $\sz S\al$ consist of all nodes $\ang{u,v}\in \sa S\al$ 
such that $\clo{\aog\al uv}\not\sq A$, where
$\aog\al uv=\pr{\pog\al uv}$, $\pog\al uv=[\sog\al uv]$, 
and
$$
\sog \al uv=\ens{\ang{s,t}\in\sa S\al}
{(u\su s\land v\su t)\lor(s\sq u\land t\sq v)}\,;
$$
\itla{+1B}
we let $\sa S{\al+1}$ be the \rit{pruning} of $\sz S\al$, 
that is, $\sa S{\al+1}$ consists of all nodes 
$\ang{u,v}\in \sz S\al$ such that 
there is an infinite branch $\ang{x,y}\in [\sz S\al]$ 
satisfying 
$u\su x$ and $v\su y$.
\een
Obviously there is a countable ordinal $\la$ such that 
$\sa S {\la+1}=\sa S\la$. 
\vom

\rit{Case 1}: $\sa S\la=\pu$. 
Then, if $x\in A=\pr P$ then by construction there exist an 
ordinal $\al<\la$ and a node $\ang{u,v}\in \sa S\al$ such 
that 
$$
x\in \aog \al uv\sq \clo{\aog \al uv}\sq A\,,
$$ 
and hence $A$ is a countable union of sets $F\sq A$ of the 
form $\clo{\aog \al uv}$, where $\al<\la$ and 
$\ang{u,v}\in \sa S\al$, closed, therefore 
\ddd\fsg compact by the above. 

Let us show how this leads to \ref{pt1} of the theorem.

First of all, quite obviously there is a certain $\is12$ 
formula $\vpi(\cdot,\cdot,\cdot)$ such that we have 
$\sa S{\al+1}=\ens{\ang{u,v}}{\vpi(\sa S\al,u,v)}$ for all 
$\al$. 
It follows by Shoenfield that the construction is absolute 
for every transitive model containing all countable 
ordinals, in particular, for $\dL$, the class of 
G\"odel constructible sets. 
Thus we can assume it from the beginning that 
\rit{we argue in\/ $\dL$}. 

Another consequence of the existence of $\vpi$ is that 
both the ordinal $\la$ and the 
sequence $\ens{\ang{\al,\sa S\al}}{\al<\la}$ are $\id13$. 
It follows 
(here we use the assumption that we argue in $\dL$) that 
each ordinal $\al<\la$ is $\id13$ and each tree 
$\sa S\al\yt\al<\la$, is $\id13$ either, as well as all 
subtrees of the form $\sog\al uv$ 
(where $\ang{u,v}\in \sa S\al$) 
and their \lap{projections} 
$\tog\al uv=\ens{u}{\sus v\,(\ang{u,v}\in \sog\al uv)}
\sq\nse$. 

On the other hand, $\clo{\aog \al uv}=[\tog\al uv]$ 
holds by construction.

To conclude, if $x\in A$ then there is a pruned $\id13$ 
tree $T\sq\nse$ 
(of the form $\tog\al uv$) 
such that $x\in [T]\sq A$ --- and $[T]$ is \ddd\fsg compact 
in this case.

It remains to note that if $T\sq\nse$ is a pruned $\id13$ 
tree and the set $[T]$ is \ddd\fsg compact then by 
Theorem~\ref{mt} (relativized version) 
there is a sequence of compact $\id11(T)$ trees $T_n$ 
such that $[T]=\bigcup_n[T_n]$. 
But each $T_n$ then is $\id13$ as so is $T$ itself. 
\vom

\rit{Case 2}: $S^{(\la)}\ne\pu$, and then $S^{(\la)}\sq S$ 
is a pruned tree. 

\ble
\lam{tt''}
If\/ $\ang{u,v}\in \sa S\la$, $u'\in\nse$, $u\su u'$, and\/ 
$\aog \la uv\cap\ibn{u'}\ne\pu$ then there is a string\/ 
$v'\in\nse$ such that\/ $v\su v'$ and\/ 
$\ang{u',v'}\in \sa S\la$.\qed
\ele

We'll define a pair $\ang{u(t),v(t)}\in \sa S\la$ 
for each $t\in\nse$, such that  
%
\ben
\itsep
\tenu{(\arabic{enumi})}
\itla{px1}
if $t\in\nse$ then $t\sq u(t)$;

\itla{px2}
if $s,t\in\nse$ and $s\sq t$ then 
$u(s)\sq u(t)$ and $v(s)\sq v(t)$;

\itla{px3}
if $t\in\nse$ and $k\ne n$ then 
$u(t\we k)$ and $u(t\we n)$ are \ddd\sq incomparable;

\itla{px4}
if $s\in\nse$ then there exists a point 
$y_s\in \clo{\aog \la {u(s)}{v(s)}}\bez A $ 
such that
any sequence of points 
$x_k\in \aog \la {u(s\we k)}{v(s\we k)}$ 
converges to $y_s$.
\een

Suppose that such a system of sets is defined. 
Then the associated map 
$f(a)=\bigcup_n u(a\res n):\bn\to A$ is $1-1$ and 
is a homeomorphism 
from $\bn$ onto its full image 
$Y=\ran f =\ens{f(a)}{a\in\bn}\sq A$.  

Let's prove that $Y$ is relatively closed in $A$. 
Consider a sequence of points $a_n\in\bn$ such that 
the corresponding sequence of points $y_n=f(a_n)\in Y$ 
converges to a point $y\in\bn$; 
we have to prove that $y\in Y$ or $y\nin A$.
If the sequence $\sis{a_n}{}$ contains a subsequence 
convergent to $b\in\bn$ then 
$\sis{y_n}{}$ converges to $f(b)\in Y$. 
So suppose that the sequence $\sis{a_n}{}$ 
has no convergent subsequences. 
Then there exist a string $s\in\nse$, 
an infinite set $K\sq\dN$, 
and for each $k\in K$ --- a number $n(k)$, 
such that $s\we k\su a_{n(k)}$. 
Then $y_{n(k)}\in \aog\la{u(s\we k)}{v(s\we k)}$ 
by construction. 
Therefore the subsequence $\sis{y_{n(k)}}{k\in\dN}$ 
converges to a point $y_s\nin A$ by \ref{px4}, 
as required.
 
Finally on the construction of sets $Y_s$. 

Suppose that a pair $\ang{u(t),v(t)}\in \sa S\la$ is defined. 
Then $\clo{\aog\la{u(t)}{v(t)}}\not\sq A$ by the choice of $\la$.
There is a sequence of pairwise different points 
$x_n\in\aog\la{u(t)}{v(t)}$ which converges to 
a point $y_s\in \clo{\aog\la{u(t)}{v(t)}}\bez A$. 
We can associate a string $u_n\in\nse$ with each $x_n$ 
such that $u(t)\su u_n\su x_n$, the strings $u_n$ are 
pairwise \ddd\sq incompatible, and $\lh{u_n}\to\infty$.
Then, by Lemma~\ref{tt''}, for each $n$ there is a 
matching string $v_n$ such that $v(t)\su v_n$ and  
$\ang{u_n,v_n}\in \sa S\la$. 
Put $u(t\we n)=u_n$ and $v(t\we n)=v_n$ for all $n$.
\vom

\epF{Theorem~\ref{pt}}


\punk{Remarks}
\las{zz}

The main results of this note can be compared with the 
following theorems of classical descriptive set theory.

\bte
[\rm Hurewicz~\cite{hur}]
\lam{hur}
If a\/ $\fs11$ set\/ $A$ in a Polish space\/ $\pX$ 
is not\/ \ddd\fsg compact then  
there is a subset\/ $Y\sq A$ homeomorphic to the 
Baire space\/ $\bn$ 
and relatively closed in\/ $A$.\qed
\ete

\bte
[\rm Saint~Raymond~\cite{sr}, see also 21.23 in \cite{dst}]
\lam{hur2}
If a\/ $\fs11$ set\/ $A$ in a Polish space\/ $\pX$ cannot 
be covered by a\/ \ddd\fsg compact set\/ $Z\sq\pX$ then 
there is a set\/ $P\sq A$,   
homeomorphic to\/ $\bn$ and closed in\/ $\pX$.\qed
\ete

Arguments in \cite{dst} show that it's sufficient to prove 
either of these theorems in the case $\pX=\bn$; then the 
results can be generalized to an arbitrary Polish space $\pX$ 
by purely topological methods. 
In the case $\pX=\bn$, Theorem~\ref{hur2} immediately follows  
from our Theorem~\ref{mt'} 
(in relativized form, \ie, 
for classes $\is11(p)$, where $p\in\bn$ is arbitrary),  
while Theorem~\ref{hur} follows from  
Theorem~\ref{pt} (relativized). 
On the other hand, Theorem~\ref{hur} also follows from  
Theorem~\ref{mt} (relativized) 
for sets $A$ in $\fd11$ (that is, Borel sets). 

Theorem~\ref{mt} implies yet another theorem, 
which combines several classical results of descriptive 
set theory by Arsenin, Kunugui, Saint~Raymond, 
She\-gol\-kov, see references in \cite{dst} or 
in \cite[\S\,4]{umnL}.

\bte
[\rm compare with Fact~\ref{25}]
\lam{tks}
Suppose that\/ $\pX,\pY$ are Polish spaces, $P\sq\pX\ti\pY$  
is a\/ $\fd11$ set, and all cross-sections\/ 
$\seq Px=\ens{y}{\ang{x,y}\in P}$ {\rm($x\in\pX$)} 
are\/ \ddd\fsg compact. 
Then\/ 
\ben
\renu
\itla{tks2}
the projection\/ $\pr P$ is a\/ $\fd11$ set$;$

\itla{tks3} 
$P$ is a countable union of\/ $\fd11$ sets with compact 
cross-sections$;$

\itla{tks1}
$P$ can be uniformized by a\/ $\fd11$ set.
\een
\ete
\bpf[a sketch for the case $\pX=\pY=\bn$]
\ref{tks2}
Assume, for the sake of simplicity, that $P\sq\bn\ti\bn$ is   
a $\id11$ set. 
The set 
$$
H=\ens{\ang{x,T}} 
{x\in\bn\land T\in\ct\land T\in\id11(x)\land [T]\sq \seq Px}
$$
is $\ip11$ by Fact~\ref{dp}. 
It follows from Theorem~\ref{mt} that if $\ang{x,y}\in P$ 
then there is a tree $T$ such that $\ang{x,T}\in H$ 
and $y\in[T]$. 
Therefore the $\ip11$ set
$$
E=\ens{\ang{x,y,T}}
{\ang{x,y}\in P\land \ang{x,T}\in H\land y\in [T]}
\sq\bn\ti\bn\ti 2^{(\nse)}
$$
satisfies $\pr_{xy}E=P$, 
that is, if $\ang{x,y}\in P$ then there is a tree $T$ such that 
$\ang{x,y,T}\in E$. 
There is a uniform $\ip11$ set $U\sq E$ which 
\rit{uniformizes} $E$, 
\ie, if $\ang{x,y}\in P$ then there is a unique $T$ such that 
$\ang{x,y,T}\in U$. 
Yet $U$ is $\is11$ as well by Fact~\ref{BQ}, 
since $\ang{x,y,T}\in U$ is 
equivalent to:
$$
\ang{x,y}\in P\land y\in [T]\land \kaz T'\in\id11(x)\:
(\ang{x,y,T'}\in U\imp T=T')\,.
$$ 
Thus the $\is11$ set 
$F=\ens{\ang{x,T}}{\sus y\:(\ang{x,y,T}\in U)}$ 
is a subset of the $\ip11$ set $H$. 
Fact~\ref{21} implies that there is a $\id11$ set $V$ such 
that $F\sq V\sq H$.
Then
$$
\ang{x,y}\in P\leqv \sus T\:
(\ang{x,T}\in V\land y\in[T]) 
$$
by definition.
Finally all cross-sections of $V$ are at most countable: 
indeed if $\ang{x,T}\in V$ then $T\in\id11(x)$ 
(since $V\sq H$). 
Note that $\pr P=\pr V$, and hence the projection $D=\pr P$ 
is $\id11$ (hence Borel) by Fact~\ref{25}.\vom 

\ref{tks3} 
It follows from Fact~\ref{25} that $V$ is equal to a  
union $V=\bigcup_nV_n$ of uniform $\id11$ sets $V_n$, 
and then each projection $D_n=\pr V_n\sq D$ is $\id11$. 
Each $V_n$ is basically the graph of a $\id11$ map  
$\tau_n:D_n\to \ct$, and 
$\seq Px=\bigcup_{x\in D_n}[\tau_n(x)]$. 
If $n\in\dN$ then we put  
$$
P_n=\ens{\ang{x,y}}
{x\in D_n\land  y\in[\tau_n(x)])}\,. 
$$
Then $P=\bigcup_nP_n$ by the above, each set $P_n$ has only  
compact cross-sections, and each $P_n$ is a $\id11$ set,  
since the sets $D_n$ and maps $\tau_n$ belong to $\id11$.\vom

\ref{tks1} 
Still by Fact~\ref{25}, the set $V$ can be uniformized by  
a uniform $\id11$ set, that is, there exists a  
$\id11$ map $\tau:D\to \ct$ 
such that $\ang{x,\tau(x)}\in V$ for all $x\in D$.
To uniformize the original set $P$, let $Q$ consist of all pairs 
$\ang{x,y}\in P$ such that $y$ is the lexicographically leftmost 
point in the compact set $[\tau(x)]$. 
Clearly $Q$ uniformizes $P$.
To check that $Q$ is $\id11$, note that 
``$y$ is a the lexicographically leftmost point in $[T]$'' 
is an arithmetic relation in the assumption that $T\in\ct$.
\epf

Similar arguments, this time based on Theorem~\ref{mt'}, also 
lead to an alternative proof of the following known result. 

\bte
[Burgess, Hillard, 35.43 in \cite{dst}] 
If\/ $P$ is a\/ $\fs11$ set in the product\/ $\pX\ti\pY$ of 
two Polish spaces\/ $\pX$, and every section\/ $\seq Px$  
is covered by a\/ \ddd\fsg compact set, then 
there is a sequence of Borel sets\/ $P_n\sq \pX\ti\pY$ 
with compact sections\/ $\seq{P_n}{x}$ such that\/    
$P\sq\bigcup_nP_n$.\qed
\ete

But at the moment it seems that no conclusive theory 
of $\fs11$ sets with \ddd\fsg compact sections 
(as opposed to those with sections covered by 
\ddd\fsg compact sets) 
is known. 
For instance what about effective decompositions of such 
sets into countable unions of definable sets with compact 
sections? 
Our Theorem~\ref{pt} can be used to show that such a 
decomposition is possible, but the decomposing sets 
with compact sections appear to be excessively complicated 
(3rd projective level by rough estimation). 
It is an interesting {\ubf problem} to improve this result 
to something more reasonable like Borel combinations of 
$\fs11$ sets.

On the other hand, it is known from \cite{sri,sti} that 
$\fs11$ sets with \ddd\fsg compact sections are not 
necessarily decomposable into countably many $\fs11$ 
sets with compact sections.

\end{document}